\documentclass{article}
\usepackage[a4paper,left=3cm,top=3.5cm]{geometry}
\usepackage{verbatim}
\usepackage[all,cmtip]{xy}

\usepackage{mathtools}

\usepackage{latexsym,bm,easymat}
\usepackage{amsmath,amscd,amssymb,amsfonts,mathdots,ntheorem}
\usepackage{indentfirst}
\usepackage[colorlinks,linkcolor=blue,citecolor=red]{hyperref}
\usepackage{latexsym}
\usepackage{longtable}
\usepackage{cases}
\usepackage{slashed}

\usepackage{shorttoc}

\begin{document}
\makeatletter
\DeclareRobustCommand\widecheck[1]{{\mathpalette\@widecheck{#1}}}
\def\@widecheck#1#2{%
    \setbox\z@\hbox{\m@th$#1#2$}%
    \setbox\tw@\hbox{\m@th$#1%
       \widehat{%
          \vrule\@width\z@\@height\ht\z@
          \vrule\@height\z@\@width\wd\z@}$}%
    \dp\tw@-\ht\z@
    \@tempdima\ht\z@ \advance\@tempdima2\ht\tw@ \divide\@tempdima\thr@@
    \setbox\tw@\hbox{%
       \raise\@tempdima\hbox{\scalebox{1}[-1]{\lower\@tempdima\box
\tw@}}}%
    {\ooalign{\box\tw@ \cr \box\z@}}}
\makeatother

\def\comp{\ensuremath\mathop{\scalebox{.6}{$\circ$}}}
\def\QEDclosed{\mbox{\rule[0pt]{1.3ex}{1.3ex}}} %
\def\QEDopen{{\setlength{\fboxsep}{0pt}\setlength{\fboxrule}{0.2pt}\fbox{\rule[0pt]{0pt}{1.3ex}\rule[0pt]{1.3ex}{0pt}}}}
\def\QED{\QEDopen} %
\def\pf{\noindent{\bf Proof}} %
\def\endpf{\hspace*{\fill}~\QED\par\endtrivlist\unskip \hfill}
\def\Dirac{\slashed{D}}
\def\Hom{\mbox{Hom}}
\def\ch{\mbox{ch}}
\def\Ch{\mbox{Ch}}
\def\Tr{\mbox{Tr}}
\def\End{\mbox{End}}
\def\Re{\mbox{Re}}
\def\cs{\mbox{Cs}}
\def\spinc{\mbox{spin}^\mathbb C}
\def\spin{\mbox{spin}}
\def\Ind{\mbox{Ind}}
\def\Hom{\mbox{Hom}}
\def\ch{\mbox{ch}}
\def\cs{\mbox{cs}}
\def\nTbe{\nabla^{T,\beta,\epsilon}}
\def\nbe{\nabla^{\beta,\epsilon}}
\def\mfa{\mathfrak{a}}
\def\mfb{\mathfrak{b}}
\def\mfc{\mathfrak{c}}
\def\mcM{\mathcal{M}}
\def\mcC{\mathcal{C}}
\def\mcB{\mathcal{B}}
\def\mcG{\mathcal{G}}
\def\mcE{\mathcal{E}}
\def\mfs{\mathfrak s}
\def\mcWs{\mathcal W_{\mathfrak s}}
\def\umcF{\underline{\mathcal F}}
\def\mce{\mathfrak e}
\def\mbT{\mathbf T}
\def\mbc{\mathbf c}
\def\spinc{$\mbox{spin}^c~$}

\newcommand{\at}{\makeatletter @\makeatother}
\renewcommand{\contentsname}{\center{Content}}
\newtheorem{defi}{Definition}[section]
\newtheorem{thm}[defi]{\textbf{Theorem}}
\newtheorem{conj}[defi]{\textbf{Conjecture}}
\newtheorem{cor}[defi]{\textbf{Corollary}}
\newtheorem{lemma}[defi]{\textbf{Lemma}}
\newtheorem{exa}[defi]{\textbf{Example}}
\newtheorem{prop}[defi]{Proposition}
\newtheorem{assump}[defi]{Assumption}
\newtheorem{hypo}[defi]{Hypothesis}

\title{Transverse $\mathcal F^T$-entropy and transverse Ricci flow for Riemannian foliations}
\date{}
\author{Dexie Lin}
\maketitle

\begin{abstract}
  In this paper, we introduce an entropy functional on Riemannian foliation, inspired by the work of , which is monotonically along the transverse Ricci flow. We relate their gradient flow, via diffeomorphism preserving the foliated structure of the manifold with Riemannian foliation, to the transverse Ricci flow. Moreover, inspired by the work of Fuquan Fang and Yuhao Zhang, we give a necessary condition for codimension 4 Riemannian foliation admitting the transverse Einstein  metric.
\end{abstract}




 \section{Introduction}

 A well known method to get a Einstein metric on Riemannian manifolds is by Ricci flow, which was originally introduced by Hamilton \cite{Hamilton}. After being used
by Perelman in proving the famous Poincar\'e conjecture,
Ricci flow became a popular topic of geometric analysis. A fundamental new discovery of Perelman's work is to prove that the
Ricci-flow evolution equation is the gradient flow of the  Perelman's $\lambda$-functional
on a Riemannian manifold \cite{Pe1}. Ricci flow was generalized to foliated Riemannian manifolds in \cite{LMR}, and
then applied to Sasaki manifolds in the paper \cite{SWZ}. Later Collins \cite{Collins} generalized Perelman's entropy functionals to Sasaki manifolds and gave the $C^0$ bound  for the
transverse scalar curvature, and a uniform $C^1$ bound for the transverse Ricci potential along the Sasaki-Ricci flow.

 On the other hand, from the viewpoint of Riemannian foliation, Sasaki manifold is a special case of Riemannian foliation, i.e. 1 dimensional foliation with taut bundle like metric. In this paper, we generalize the transverse $\lambda$-functional to  general Riemannian foliations. More precisely, let $(M,F)$ be a closed manifold $M$ with a Riemannian foliation $F$,
   we introduce a new functional for a bundle like metric $g$ and a smooth basic function $f$ as follows,
   \[{\mathcal F}^T(g,f)=\int_M(Scal^T+|\nabla f|^2+ |\kappa_b|^2+2\tau_bf)e^{-f},\]
   where $Scal^T$ is the transverse scalar curvature, $\tau_b$ is the basic mean curvature field and $\kappa_b$ is dual to $\tau_b$ via the bundle-like metric. We refer the readers to Section 2   for the relevant definitions.

\noindent
    In this paper, we fix the leafwise metric $g_F$ and deform the transverse metric $g^T$.
    Similar to the manifold case, we define the number $$\lambda^T(g)=\min\{{\mathcal F}^T(g,f)|\int_M e^{-f}dVol=1\}.$$
   We have the following proposition about this functional. Set
   \[\bar\lambda^T(g)=\lambda^T(g)Vol(g)^{2/m},\]
   where $m$ is the codimension of  $F$. This is invariant under the scaling of the transverse metric.
   \begin{prop}\label{prop-main}
     Let $(M,F)$ be a closed manifold with a Riemannian foliation $F$. Let $g(0)$ be the initial bundle like metric. Suppose that we fix the leafwise metric $g_F=g\big|_F$ and deform the transverse part along the flow
     \begin{equation}
       \frac{d}{dt}g^T(t)=-2Ric^T,\label{eqn-transverse-Ricci-flow}
     \end{equation}
     where $g^T$ denotes the transverse part of the metric $g$. Then, the above functional $\lambda^T$ is weakly non-decreasing along this flow.
   \end{prop}

\noindent
   It is known that for manifolds with  codimension 4 Riemannian foliations, under a certain topological  condition(see \cite{KLW}), one can define the baisc Seiberg-Witten equations.
   We will give a necessary condition for Riemannian foliations  admitting a transverse Einstein metric.

  \begin{thm}\label{thm-main-1}
      Let $(M,F)$ be an oriented closed codimension $4$ foliation. Suppose  that $(M,F)$ admits a bundle-like metric $g$ and a transverse $spin^c$ structure  $\mathfrak s$. If there is an irreducible solution $(A,\psi)$ to the basic Seiberg-Witten equations with respect to $g$ and $\mathfrak s$, then  \[\bar\lambda^T(g)\leq -\sqrt{\|F^+_A\|^2}.\]
      If $[c^+_b(\mathfrak s)]\neq0$, then  the equality holds if  and only if $g$ is a transverse K\"ahler metric and is taut bundle like metric, i.e. each leaf is minimal with respect to this metric.
   \end{thm}

   \noindent
   We refer the readers to Section 3  for the relevant definitions. Combining with  the similar arguments of \cite{FZ}, we have the following theorem.

   \begin{thm}\label{thm-main-2}
      Let $(M,F, \mathfrak s)$ be the triple of the above theorem. Suppose that there is a bundle like  metric $g$ admitting  an irreducible solution to the basic Seiberg-Witten equations.  If $[c^2_b(\mathfrak s)]>0$, then for any bundle like metric $g$, we have that
      \[
      \bar\lambda^T(g)\leq-\sqrt{32\pi^2 c^2_b(\mathfrak s)},\]
      the equality holds if and only if $g$ is transverse K\"ahler Einstein and taut bundle like.
   \end{thm}

\noindent
 The organization of this paper is as follows: in Section 2, we review the some necessary results for manifold with Riemannian foliations, and give the definition of basic entropy functional; in Section 3, we review the work of basic Seiberg-Witten equations and give the proofs of our main theorem.

 \vspace{3mm}

\noindent
 {\bf Acknowledgement}:
 The author would express the special appreciation  to Ken Richardson for the discussing about the transverse diffeomorphism.

 \section{Riemannian foliations and transverse Ricci flow}

 In this section, we first review some results of the previous work on Riemannian foliations and give the definition of transverse $\mathcal F^T$ entropy functional for manifolds with a Riemannian foliation.

 \subsection{Geometry of Foliations}

  In this subsection, we let $M$ be a closed  $n$ dimensional manifold with dimension $p$ foliation $F$. We  denote codimension of this foliation by $m=n-p$. For more details of this subsection, we give a reference \cite{Ton2}.
 \begin{defi}
   A Riemannian metric $g^T$ on $Q$ is said to be bundle-like, if
   \[L_Xg^T\equiv0,\]
   for any $X\in \Gamma(F)$, where $Q=TM/F$. We say $(M,F)$ is a Riemannian foliation, if $Q=TM/F$ admits a bundle-like metric.
 \end{defi}
 Given a metric $g$ on $TM$, $Q$ can be identified with the orthogonal complement  $F^\perp$ induced by $g$. In turn, $Q$ inherits a metric $g_{F^\perp}$. We have the following equivalence,
 \[\mbox{a metric }g\mbox{ of }TM\mbox{ corresponds a triple }(g_F,\pi_F, g^T),\]
 where $g_F=g|_{F}$ and $\pi_F$ is the projection $TM\to TF$.
   Let $M$ be a manifold with foliation $F$.  A Riemannian metric $g$ on $TM$ is said \emph{bundle-like}, if the induced metric $g^T$ is bundle-like.

\noindent
 By the work of Reinhart $\cite{Reinhart}$, it is known  that the bundle-like metric can be locally written as $g=\sum_{i,j}g_{ij}(x,y)\omega^i\otimes \omega^j+\sum_{k,l}g_{k,l}(y)dy^k\otimes dy^l$, where $(x,y)$ is in the foliated chart of $M$ and $\omega^i=dx^i+a^i_\alpha(x,y) dy^\alpha$.
In this paper, we assume that $(M,F)$ is a Riemannian foliation.
Let $\pi$ be the canonical projection $TM\to Q$.  We define a connection $\nabla^{T}$ on $Q$, by
$$\nabla^{T}_Xs:=\begin{cases}
   \pi([X,Z_s])& X\in \Gamma(F),\\
   \pi(\nabla_X Z_s)& X\in \Gamma(F^\perp),
 \end{cases}$$
 for any section $s\in\Gamma(Q)$,
 where $Z_s\in \Gamma(TM)$ is a lift of $s$, i.e. $\pi(Z_s)=s$.  One can verify that it is  torsion free and metric-compatible, whose leafwise restriction coincides with the Bott-connection. We set $R^T$ as the curvature of this connection.
 Parallel to the Riemannian manifold, w the transverse Ricci curvature and scalar curvature are defined by
 $Ric^{T}(Y)=\sum^{q}_{i=1}R^{T}(Y,e_i)e_i,~
 Scal^{T}=\sum^{q}_{i=1}g^T(Ric^{T}(e_i),e_i),$
 where $\{e_i\}$ is a local  orthonormal frame of $Q$.  $R^T$
satisfies the   condition
$\iota_XR^T=0$,
for any field $X\in\Gamma(F)$.
We recall some notions.

\begin{defi}
  We say a  Riemannian  foliation $(M,F)$ admits a transverse Einstein metric, if there is a constant $c\in\mathbb R$ such that $$Ric^T=c g^T,$$ for some bundle-like metric $g$.
\end{defi}
  We define the basic forms as follows:
 \[\Omega^r_b(M)=\{\omega\in\Omega^r(M)\big|~\iota_X(\omega)=0,~L_X(\omega)=0,\forall X\in\Gamma(F)\}.\]
We set $d_b$ as the restriction of $d$ to the basic forms, the complex $d_b:\Omega^r_b(M)\to\Omega^{r+1}_b(M)$ is a subcomplex of the deRham complex, whose cohomology is called  basic cohomology, and  denoted by $H^r_b(M)$. It is known that $H^1_b(M)\subset H^1(M)$.
 We denote by $b^r_b=\dim H^r_b(M)$.
\begin{defi}
  The mean curvature vector field is defined by $\tau=\sum^{\dim F}_{i=1}\pi(\nabla_{\xi_i}\xi_i)\Gamma(Q)$, where $\{\xi_i\}$ is a local orthonormal basis of $F$. Let $\kappa\in \Gamma(Q^*)$ be the dual to $\tau$ via the metric  $g^T$.
\end{defi}

\noindent
 By the work of Alvarez L\`opez \cite{AL}, we have the following $L^2$ orthogonal decomposition  for the forms on $M$,
 \[\Omega(M)=\Omega_b(M)\oplus\Omega_0(M).\]
 One has the basic Hodge-star operator,
 \[\bar*:\bigwedge^rQ^*\to \bigwedge^{q-r}Q^*.\]
 Choosing  a local orthonormal basis $\{e_i\}_{1\leq i\leq p}$ of $F$, we define the character form of the foliation $\chi_F$ by
 \[\chi_F(Y_1,\cdots, Y_p)=\det(g_F(e_i,Y_j))_{1\leq i,j\leq p},\]
 for any section $Y_1,\cdots, Y_p\in \Gamma(TM)$.
 The basic Hodge-star operator is  related to the usual Hodge-star operator by the formula $\bar*\alpha=(-1)^{(q-r)\dim(F)}*(\alpha\wedge\chi_F)$.
 We have  the volume density formula,
 $dvol_M=dvol_Q\wedge\chi_F$.
 For  any section $\alpha\in \Gamma(\bigwedge^rQ^*)$, we define its $L^2$ norm by
 \[\|\alpha\|^2_{L^2}=\int_M\alpha\wedge\bar*\alpha\wedge\chi_F.\]
 For a bundle-like metric, it is clear that
 \[\bar*:\Omega^r_b(M)\to\Omega^{q-r}_b(M).\]

\begin{prop}[Rummler \cite{Rummler}]\label{formula-Rummler}
For any metric $g$ on $TM$,
  we get
  \[d\chi_F=-\kappa\wedge\chi_F+\phi_0,\]
  where $\phi_0$ belongs to $F^2\Omega^p=\{\omega\in\Omega^p(M)\big|\iota_{X_1}\cdots\iota_{X_p}\omega=0,\mbox{ for any }X_1,\cdots,X_p\in \Gamma(F)\}$.
\end{prop}

\noindent
Note that the above formula implies that the mean curvature form $\kappa$ is invariant under the deformation of the metric $g^T$ on the transverse bundle.
By the decomposition, we have the decomposition  $\kappa=\kappa_b+\kappa_0$, where $(\kappa_0,\omega_b)_{L^2}=0$ for any basic one form $\omega_b$.
We call $\kappa_b$ the basic mean curvature form.  It is known that
 $d\kappa_b=0$,
 and the cohomology class $[\kappa_b]$ is independent on any bundle-like metric \cite{AL}.
It is known that the space of bundle-like metrics on $(M, F)$ is infinite
dimensional \cite{SC}. It can be shown that any bundle-like metric can be deformed in
the leaf directions leaving the transverse part untouched in such a way that the
mean curvature form becomes basic \cite{D}.

 \begin{prop}[c.f. {\cite[Theorem 7.18]{Ton2}}]
   The  $L^2$-formal adjoint   of $d_b$ is $\delta_b=(-1)^{m(*+1)+1}\bar*(d_b-\kappa_b\wedge)\bar*$.
 \end{prop}

\noindent
We define the basic Laplacian operator $\Delta_b=d_b\delta_b+\delta_bd_b$.





\begin{defi}
  We say a foliation is \emph{taut}, if  there is a metric on $M$ such that $\kappa=0$, i.e. all leaves are minimal submanifolds.
\end{defi}

\noindent
For a fixed Riemannian foliation $F$, the taut condition has a topological obstruction.

\begin{prop}[Alvarez L\`opez \cite{AL}]
  Let $F$ be a  Riemannian foliation on a closed manifold $M$. Then, $F$ is taut if and only if the class $[\kappa_b]$ is trivial. Furthermore, when $F$ is transversely oriented  the foliation is taut if and only if $H^m_b(M)\neq0$.

\end{prop}

\subsection{Transverse Ricci flow}
 Similar to Hamilton's arguments on the Ricci flow, the short existence of  a transverse Ricci flow\[\frac{d}{dt}g^T(t)=-2Ric^T,\]
 was mentioned in \cite{LMR} for a  Riemannian  foliation.
\noindent
In general, 
 for any bundle-like metric $g=g_F\oplus g^T$, by using the Einstein convention  we have that
\[R^{T,l}_{ijk}=\partial_i\Gamma^{T,l}_{jk}-\partial_j\Gamma^{T,l}_{ik}
+\Gamma^{T,p}_{jk}\Gamma^{T,l}_{ip}-\Gamma^{T,p}_{ik}\Gamma^{T,l}_{jp},\]
where $\Gamma^{T}$ is the (transverse)Christoffel symbol for the transverse Levi-Civita connection.
since the transverse curvature $R^T$ also has the first and second Bianchi identities.
For a variation $\frac{\partial}{\partial t}g^{T}=v$, we have that
\[d_vRic^T_{ij}=\frac12(\nabla^T_l(\nabla^T_i)v_{jl}+\nabla^T_jv_{il}
-\nabla^T_lv_{ij})-\frac12 \nabla^T_i\nabla^T_jV,\]
where $V=tr(v)=\sum_{ij}g^{T,ij}v_{ij}$ and $d_v$ denotes the variation with respect to $v$.
Hence, one obtains that
\[d_vScal^T=div(div(v))-\Delta'_bV-(v,Ric),\]
where $div(div(v))=g^{T,ij}\nabla^T_i\nabla^T_jv_{ij}$ and $\Delta'_b=-\Delta_b+\nabla^T_{\tau_b}$
We define  the  ${\mathcal F}^T $ entropy functional by
\begin{equation}
  {\mathcal F}^T(g,f)=\int_M(Scal^T_g+|\nabla f|^2+|\kappa_b|^2+2\tau_bf)e^{-f}\label{eqn-mu-F-functional}
\end{equation}for a bundle-like metric $g$ and a basic function $f$ on a foliated manifold $(M,F)$. Under the constrain $\{\int_M e^{-f}dVol=1\}$,
we define the number $\lambda^T(g)=\min\{{\mathcal F}^T(g,f)|\int_M e^{-f}dVol=1\}$. By setting $\Phi=e^{-f/2}$,  $\lambda^T(g)$ is the first eigenvalue of the operator
\[4\Delta_b+(Scal^T+|\kappa_b|^2-2\delta_b\kappa_b).\]
  Before proceeding, we need the following  lemma to assist the later calculation.

\begin{lemma}
  Under the above notations, for any  $v\in \Gamma_b(Q^*\odot Q^*)$ and any $Y\in \Gamma_b(Q^*)$, we have that
  \[\int_M(div(v),Y))=\int_M(v,Y\odot\kappa_b)-\int_M(v,\nabla^TY),\]
  where $Q^*\odot Q^*$ denotes the symmetric tensor of $Q^*\otimes Q^*$, $div(v)=\sum^m_{i=1}\iota_{e_i}\nabla^T_{e_i}v$ and $\kappa_b\odot Y=\frac12(\kappa_b\otimes Y+Y\otimes \kappa_b)$.
\end{lemma}
\begin{pf}
  By straightforward calculation, we have that
  \begin{eqnarray*}
    \int_M(div(v),Y))&=&\int_M(\sum^q_{i=1}\iota_{e_i}\nabla^T_iv,Y)dvol\\
    &=&\int_M(\sum^m_{i=1}\nabla^T_{e_i}v,e^i\odot Y)dvol\\
    &=&\int_M\sum^m_{i=1}e_i(v,e^i\odot Y)dvol-\int_M\sum^m_{i=1}(v,e^i\odot \nabla^T_{e_i}Y)dvol.
  \end{eqnarray*}
  We define $X\in \Omega^1_b(M)$,  such that $(X,Y)=(v,X\odot Y)$, for any $Y\in\Omega^1_b$. Then, it holds that
  \[\int_M div(X)=\int_M\sum^m_{i=1}(\nabla^T_{e_i}X,e^i)=\int_M\sum^m_{i=1}
  e_i(X,e^i)=\int_M\sum^m_{i=1}e_i(v,e^i\odot X)=\int_{M}(v,\kappa_b\odot Y),\]
  where the last equality comes from \cite[Theorem 5.24]{Ton2}.
  Hence, one deduces that
  \[\int_M (div(v),Y)=\int_M(v,\kappa_b\odot Y)-\int_M(v,\nabla^T Y).\]
\end{pf}

Before proceeding, we need a full variation of $\mathcal F^T$.
\begin{lemma}
  Under the above conditions, we have that the formula below
\begin{eqnarray*}
    &&d_{(v,h)}\mathcal F^T(g,f)\\
    &=&\int_M-(v,Ric^T_g+Hess_b(f)+\nabla^T\kappa_b)(e^{-f}dvol)\\
    &&+
    \int(\frac{V}{2}-h)\left((Scal^T_g-2\Delta_bf-|\nabla^Tf|+|\kappa_b|^2-2\delta_b\kappa_b))\right)e^{-f}dvol,
  \end{eqnarray*}
where $v$ is a basic quadratic tensor on $(M,F)$, $h$ is a basic function and $d_{(v,h)}$ denotes the variation of $\mathcal F^T$ with respect to $(v,h)$.
\end{lemma}
\begin{pf}
 Using the previous lemma, one deduces that
\begin{eqnarray*}
  &&\int_Mdiv(div(v))e^{-f}\\
  &=&\int_M(div(v),\kappa_b)e^{-f}-\int_M
  (\nabla^Te^{-f},div(v))\\
  &=&\int_M(div(v),\kappa_b)e^{-f}-\int_M(v,\nabla^T e^{-f}\odot\kappa_b)+
  \int_M(v,Hess(e^{-f}))\\
  &=&\int_M(v,Hess(e^{-f}))+\int_M(v,\kappa_b\odot\kappa_b)e^{-f}-2\int_M
  (v,\nabla^T(e^{-f})\odot\kappa_b)-\int_M(v,\nabla^T\kappa_b)e^{-f}.
\end{eqnarray*}
We rewrite the variation as two parts
\[d_{(v,h)}{\mathcal F}^T(g,f)=\int_Md_{(v,h)}(Scal^T+|\nabla f|^2+ |\kappa_b|^2+2\tau_bf)e^{-f}dvol+
\int_M(Scal^T+|\nabla f|^2+ |\kappa_b|^2+2\tau_bf)d_{(v,h)}(e^{-f}dvol).
\]
By the identity $\tau_bf=(\kappa_b,df)$,  we have
\[2\int_Md_{(v,h)}(\tau_bf)e^{-f}dvol=-2\int_M (v,df\otimes\kappa_b)e^{-f}dvol+\int_M\tau_b Ve^{-f}dvol=
2\int_M (v,\nabla^Te^{-f}\otimes\kappa_b)dvol+\int_M\tau_b Ve^{-f}dvol,\]
and
\[\int_M d_{(v,h)}(|\kappa_b|^2)e^{-f}dvol
=\int_M-(v,\kappa_b\otimes\kappa_b)e^{-f}dvol.\]
Similarly, we deduce  the formulas
\begin{eqnarray*}
  \int_M(-\Delta'_bV)e^{-f}=
  \int_M\Delta_bVe^{-f}-\tau_bVe^{-f}=\int_M(\Delta_be^{-f})V-\tau_bVe^{-f}
\end{eqnarray*}
and
\begin{eqnarray*}
  &&\int_Md_{(v,h)}(|\nabla^Tf|^2)e^{-f}dvol\\
  &=&\int_M-v(\nabla^Tf,\nabla^Tf)e^{-f}+2(\nabla^Tf,\nabla^T\frac{\partial f}{\partial t})e^{-f}\\
  &=&\int_M-v(\nabla^Tf,\nabla^Tf)e^{-f}+(\nabla^Tf,\nabla^TV)e^{-f}\\
  &=&\int_M-v(\nabla^Tf,\nabla^Tf)e^{-f}-(de^{-f},dV)\\
  &=&\int_M-v(\nabla^Tf,\nabla^Tf)e^{-f}-\delta_bd(e^{-f})V=\int_M-v(\nabla^Tf,\nabla^Tf)e^{-f}
  -(\Delta_be^{-f})V.
\end{eqnarray*}
Combining the above calculations, one has
  \begin{eqnarray*}
    &&\int_Md_{(v,h)}(Scal^T_g+|\nabla^Tf|^2+|\kappa_b|^2+2\tau_bf)(e^{-f}dvol)\\
    &=&\int_M-(v,Ric^T_g+Hess_b(f)+\nabla^T\kappa_b)(e^{-f}dvol)\\
    &&+\int_M(-\Delta'_bV)e^{-f}-2\int_M(\Delta_be^{-f})hdvol+2\int_M\tau_bh(e^{-f}dvol)\\
    &=&\int_M-(v,Ric^T_g+Hess_b(f)+\nabla^T\kappa_b)(e^{-f}dvol)\\
    &&+2\int_M(\Delta_be^{-f})(\frac{V}{2}-h)dvol-2\int_M\tau_b(\frac{V}{2}-h)(e^{-f}dvol).
  \end{eqnarray*}
  For the measure part, it holds that
  \begin{eqnarray*}
    \int_M(Scal^T_g+|\nabla^Tf|^2+|\kappa_b|^2+2\tau_bf)d_{(v,h)}(e^{-f}dvol)
    =\int_M(Scal^T_g+|\nabla^Tf|^2+|\kappa_b|^2+2\tau_bf)(\frac{V}{2}-h)(e^{-f}dvol).
  \end{eqnarray*}
  By the elementary calculation, one has
  \begin{eqnarray*}
    \int_M\tau_b(\frac{V}{2}-h)(e^{-f}dvol)&=&\int_M(d(\frac{V}{2}-h),e^{-f} \kappa_b)dvol
    =\int_M
  (\frac{V}{2}-h)\delta_b(e^{-f}\kappa_b)dvol\\
  &=&\int_M(\frac{V}{2}-h)\delta_b(\kappa_b)e^{-f}dvol+\int_M(\frac{V}{2}-h)(\tau_bf)e^{-f}dvol.
  \end{eqnarray*}
  By the identity $\Delta_b(e^{-f})=\Delta'_b(e^{-f})+\tau_b(e^{-f})=(-\Delta_bf)e^{-f}-|\nabla^Tf|^2e^{-f}$, we finish the proof by the following formula:

  \begin{eqnarray*}
    &&d_{(v,h)}\mathcal F^T(g,f)\\
    &=&\int_M-(v,Ric^T_g+Hess_b(f)+\nabla^T\kappa_b)(e^{-f}dvol)\\
    &&+
    \int(\frac{V}{2}-h)\left((Scal^T_g-\Delta_bf+|\kappa_b|^2-2\delta_b\kappa_b)e^{-f}+\Delta_b(e^{-f})\right)dvol\\
    &=&\int_M-(v,Ric^T_g+Hess_b(f)+\nabla^T\kappa_b)(e^{-f}dvol)\\
    &&+
    \int(\frac{V}{2}-h)\left((Scal^T_g-2\Delta_bf-|\nabla^Tf|+|\kappa_b|^2-2\delta_b\kappa_b))\right)e^{-f}dvol.
  \end{eqnarray*}


\end{pf}

When the measure $e^{-f}dvol$ is fixed,  the following holds.

\begin{cor}
Let $v$ be a basic quadratic tensor on $(M,F)$ and $h$ be a basic function satisfying $tr(v)/2=h$. Then,
 we have
\begin{eqnarray*}
  d_{v}\mathcal F^T=-\int_M(v,Ric^T+Hess_b(f)+\nabla^T\kappa_b)e^{-f}dvol.
\end{eqnarray*}

\end{cor}

\noindent
It is clear that $\mathcal F^T$ is monotonically non-decreasing under the following flow,
\begin{equation}
  \begin{cases}
  \frac{\partial }{\partial t}g^T=-2(Ric^T+Hess_b(f)+
  \nabla^T\kappa_b),\\
 \frac{\partial }{\partial t}f=-Scal^T-\Delta'_bf-\delta^T\kappa_b,
\end{cases},\label{eqn-Ricci-heat-flow-1}
\end{equation}
where $(\delta^T\kappa_b)(g)=\sum_\alpha(\nabla^g_{v_\alpha}\kappa_b)(v_\alpha)$ for a local orthonormal basis $\{v_\alpha\}$ on $Q$ and the Levi-Civita connection $\nabla^g$ with respect to the metric $g$.
Let $u$ be a basic vector field over $(M,F)$, i.e. its dual via $g $ is  a basic one-form, and let $\rho(t)$ be the one-parameter family of diffeomorphism  generated by $u$. We have the lemma below.

\begin{lemma}\label{lemma-transverse-differmorphism}
  Under the above assumption, we have that
  \begin{itemize}
    \item[(i)] $\rho(t)_*F=F$;
    \item[(ii)] $\rho(t)^*g$ is a bundle-like metric,
  \end{itemize}
  for each $t$ in a closed interval $[0,a]$.
\end{lemma}
\begin{pf}
  Fix a point $o\in M$ and choose a local coordinates$\{x_1,\cdots,x_p,y_1,\cdots,y_m\}$ in an open neighborhood $U$ of $o$ with respect to $g(0)$.
  Choose a  small enough $t$, such that $\rho(s)(x_1,\cdots,x_m,0...0)$ remains in the coordinates patch for $|x^2_i|<\epsilon,~i=1...,m$ for some $\epsilon$. Since $V $ is independent on the leaf. For each curve of $\rho(s)(x_1,...,x_m,0...,0)$, its vector locates on the same tangent space of the leaf through $o$. The result for general $t$ is similar to \cite[Lemma 4.1]{Collins}.

  For the statement $(ii)$, we recall the notion of flat coordinate chart and adapted frame, c.f. \cite[Definition 3.1]{Y}. Choose a flat coordinate chart(see \cite[Section 2]{Y}) $U(x^i,y^\alpha)$, we have the transformation under the diffeomorphism $\rho(t)$ for small enough $t$ is still   a flat coordinate chart. Let $\{v_1,...,v_p,y_1,...y_m\}$ be an adapted frame with respect to $U(x^i,y^\alpha)$.
  Hence $$\{\rho(t)^{-1}_*v_1,..\rho(t)^{-1}_*v_p,
  \rho(t)^{-1}_*y_1,...\rho(t)^{-1}_*y_m\}$$ is an adapted frame with respect to $\rho(t)U(x^i,y^\alpha)$. Let $\nabla^g$ be the Levi-Civita connection associated with  the metric $g$.  By \cite[Theorem 3.1]{Y}, it enough to check that \[(\rho^*g)(\nabla^{\rho^*g}_{\rho^{-1}_*y_\alpha}\rho^{-1}_*v_i,\rho_*y_\beta)
  +
  (\rho^*g)(\nabla^{\rho^*g}_{\rho^{-1}_*y_\beta}\rho^{-1}_*v_i,\rho^{-1}_*y_\alpha)
  =0\]
  Since $\rho^*(g(u,w))=(\rho^*g)(\rho^{-1}_*u,\rho^{-1}_*w)$ 
  and $\nabla^{\rho^*g}_wu=\rho^{-1}_*(\nabla^g_{\rho_*w}\rho_*v)$ for any two vector fields $u,~w$ on $M$, we have that the pull-back metric is still bundle like. The remaining statement of $(ii)$ is a application of statement $(i)$ and the definition of bundle-like metric.
\end{pf}

\noindent
Inheriting the notations of \cite{Collins}, we define
\[\mathfrak Met^T(M,F)=\{g\in \mathfrak Met(M)|~ g \mbox{ is a bundle like metric.}\},\]
and define the transverse diffeomorphism group $\mathfrak Diff^T(M,F)$ by
\[\mathfrak Diff^T(M,F)=\{\rho\in \mathfrak Diff(M)|~\rho_*F=F,~\rho^*(g^T)=(\rho^*g)_Q \}.\]

\begin{prop}\label{prop-F-functional-invariant}
  The functional $\mathcal F^T$ is invariant under the action of $\mathfrak Diff^T(M,F)$.
\end{prop}

\begin{pf}
 Recall that by the Koszul-formula \cite[Theorem 5.9]{Ton2},  for any section $Y\in\Gamma(TM)$, we  write the transverse connection as
 \[g^T(\nabla^T_Ys,s')=\frac12(Yg^T(s,s')+
 Z_sg^T(\pi(Y),s')-Z_{s'}g^T(\pi(Y),s)\]\[+g^T(\pi([Y,Z_s]),s')+
 g^T(\pi([Z_{s'},Y]),s)-g^T(\pi([Z_s,Z_{s'}]),Y)),\]
 where $Z_s$ and $Z_{s'}$ are the lifts of $s$ and $s'$ respectively. Thus, the connection  $\nabla^T$ and the  transverse  scalar curvature are  uniquely determined by $g^T$.
 Since the diffeomorphism $\rho$ preserves the foliation, we have that $\rho^*(g_F)=(\rho^*g)_F$, 
 hence the scalar curvature along each leaf has the property $\rho^*Scal^F({g_F})=Scal^F(\rho^*(g_F) )$.
 Let $\nabla^g$ be the Levi-Civita connection with respect to the metric $g$.
 By the formula of \cite[Corollary 2.5.19]{BG}, we have that
 \begin{eqnarray*}
   \rho^*Scal^T(g)-Scal^T(\rho^*g)&=&\rho^*|A(g)|^2-|A(\rho^*g)|^2-
   \rho^*|T(g)|^2+|T(\rho^*g)|^2\\
   &&+\rho^*|\kappa_b(g)|^2-|\kappa_b(\rho^*g)|^2+
   2\rho^*((\delta^T\kappa_b)(g))-2(\delta^T\kappa_b(\rho^*g)).
 \end{eqnarray*}
 where the terms  $|A(g)|^2$ and $|T(g)|^2$ are defined by
   $|A(g)|^2=\sum_{i,\alpha} g^T(\nabla^g_{e_i}v_\alpha,\nabla^g_{e_i}v_\alpha)$ and
 $|T(g)|^2=\sum_{i,\alpha} g_F(\nabla^g_{v_\alpha }e_i,\nabla^g_{v_\alpha }e_i)$
  for a local orthonormal basis $\{e_i\}$ on the foliation $F$ and a local orthonormal basis $v_\alpha$ on $Q$ with respect to the bundle-like metric $g$.
 By the identity $${\nabla^{\rho^*g}}_{u} w=\rho_{*}^{-1}\left(\nabla^g_{\rho_* u} \rho_{*} w\right),$$
   where $u$ and $w$ are any vector fields on $M$,
one has the following formulas:
 $$\rho^{*}(g^{T}(\nabla^g_{e_{i}} v_\alpha, \nabla^g_{e_{i}} v_\alpha))=
 (\rho^{*} g^{T})(\rho_{*}^{-1} \nabla^g_{e_{i}}v_\alpha, \rho_{*}^{-1} \nabla^g_{e_{i}}v_\alpha)$$
  and
  $$\rho^{*}\left(g_{F}\left(\nabla^g_{v_\alpha} e_{i}, \nabla^g_{v_\alpha} e_{i}\right)\right)=(\rho^{*} g_{F})\left(\rho_{*}^{-1} \nabla^g_{v_\alpha} e_{i}, \rho_{*}^{-1} \nabla^g_{v_\alpha} e_{i}\right).$$
   Combining with the fact that the mean curvature form is uniquely  determined by the metric $g_F$ on the foliation and it is invariant under the transverse diffeomorphism, we get that $\rho^*Scal^T(g)=Scal^T(\rho^*g)$.
 For the remaining terms, one can also verify that they are also invariant under $\rho$.
\end{pf}

\begin{prop}\label{prop-monotonical-nondecreasing}
  The basic eigenvalue $\lambda^T(g(t))$ of  $\mathcal F^T$ functional is s a weakly increasing function along the variant transverse Ricci flow \eqref{eqn-transverse-Ricci-flow}.
\end{prop}

\begin{pf}
  The idea is similar to the Ricci flow on the manifold case. We consider the a family of transverse diffeomorphism
  \[\rho(t):M\to M,\]
  such that
  \[\frac{\partial\rho}{\partial t}=-(\nabla^T_{g(t)}f(t)+\tau_b),~\rho(0)=Id,\]
  then the above flow \eqref{eqn-Ricci-heat-flow-1} becomes
  \begin{equation}
  \begin{cases}
  \frac{\partial }{\partial t}g^T=-2Ric^T,\\
 \frac{\partial }{\partial t}f
  =-Scal^T-\Delta'_bf+|\nabla^Tf|^2+\tau_bf-\delta^T\kappa_b.
 \end{cases}
 \label{eqn-Ricci-heat-flow-2}
  \end{equation}

  By setting $u=e^{-f}$, the second equation of \eqref{eqn-Ricci-heat-flow-2} becomes
  \begin{equation}
    \frac{\partial u}{\partial t}=\Delta_bu+(Scal^T+\delta^T\kappa_b)u.\label{eqn-backward-heat-flow}
  \end{equation}
The equations \eqref{eqn-Ricci-heat-flow-2} will be applied by the
  transverse Ricci flow $\frac{\partial }{\partial t}g^T=-2Ric^T$, on some compact interval $[0,T]$, then by reversing the time $t\mapsto T-t$ to the
  above equation \eqref{eqn-backward-heat-flow}, one has the transverse heat flow on closed manifold with Riemannian foliation, which can be always solved on the compact interval $[0,T]$ with initial function $\bar u(0)=u(T)$. Combining with the above Proposition \ref{prop-F-functional-invariant} and Lemma \ref{lemma-transverse-differmorphism}, we finish the proof.
\end{pf}

\noindent
We put $\bar\lambda^T(g)=\lambda^T(g^T)Vol(g^T)^{2/m}$, which is transversely scale-invariant. Running the similar arguments of \cite[Claim  2.3]{Pe1}, one has the the proposition below.

\begin{prop}
  Under the taut condition $\kappa_b=0$,
   for a  fixed   leafwise metric $g_F$ we have that  if $\bar\lambda(g(t))\leq0$ for some $t$  along a flow of \eqref{eqn-transverse-Ricci-flow}, then $\frac{d}{dt}\bar\lambda(g(t))\geq0$.
\end{prop}
\noindent
{\bf Remark}:
In the general case, i.e. $\kappa_b\neq0$, the author does not know whether the proposition above holds or not.


\section{Proof of the main theorem}
Before giving the proof,   we review  the notions of  foliated bundle, basic connection, basic Dirac operator and basic Seiberg-Witten equations.

\begin{defi}
   A principal bundle $P\to M$ is called foliated, if it is equipped with a lifted foliation $F_P$ invariant under the structure group action, such that it is transversal to  the tangent space to the fiber and $F_P$ projects isomorphically onto $F$. We say a vector bundle $E\to M$ is foliated, if its principal bundle $P_E$ is foliated.
 \end{defi}

 \begin{defi}
   A connection $\omega$ of the foliated principal bundle $P$ is called adapted, if the horizontal distribution associated to this connection   contains the foliation $F_P$. A covariant connection on a foliated vector bundle is called adapted, if its associated connection on the principal bundle is. We say an adapted connection $\omega$ is called basic, if it is a Lie algebra valued basic form. Similarly, an adapted covariant connection is called basic, if its principal connection is.
 \end{defi}
\noindent
 For a foliated vector bundle $E$ over $M$, by the definition, it is clear that
 for any two connections $\nabla^1$ and $\nabla^2$ adapted to this foliated vector bundle, it holds that
 \[\nabla^1_Vs=\nabla^2_Vs,\]
 for all $s\in\Gamma(M,E)$ and $V\in\Gamma(F)$. We define the basic sections by
 \[\Gamma_b(E)=\{s\in \Gamma(E)\big|~ \nabla_Xs\equiv0,~\mbox{for all }X\in \Gamma(F)\},\]
 where $\nabla$ is an adapted connection.
 \begin{defi}
   A transverse Clifford module $E$ is a complex vector bundle over $M$ equipped with a hermitian metric  satisfying the following properties.
 \begin{enumerate}
   \item  $E$ is a  bundle of $Cl(Q)$-modules, and the Clifford action $Cl(Q)$ on $E$ is skew-symmetry, i.e.  \[(s\cdot\psi_1,\psi_2)+(\psi_1,s\cdot\psi_2)=0,\]
 for any $s\in\Gamma(Q)$ and $\psi_1,\psi_2\in\Gamma(E)$;
   \item $E$ admits a basic metric-compatible connection, and this connection is compatible with the Clifford action.
 \end{enumerate}
 \end{defi}
We say $(M,F)$ admits a transverse \spinc structure \cite{PR}, if $Q$ is \spinc and associated the spinor bundle  $S$ is a  transverse Clifford module over $(M,F)$. Moreover, $S$ induces a foliated complex line bundle $\det(S)$. 
 \begin{defi}
   Fixing a basic connection $\nabla^E$, we define the Dirac operator $\Dirac_{tr}$ by $\Dirac_{tr}=\sum^q_{i=1}e_i\cdot\nabla^E_{e_i}$ action on $\Gamma(E)$, where $\{e_i\}$ is a local orthonormal basis of $Q$.
\end{defi}
 Note that $\Dirac_{tr}$ is not formally self adjoint in general, whose adjoint operator is $\Dirac^*_{tr}=\Dirac_{tr}-\tau_b$.  The operator $\Dirac_b=\Dirac_{tr}-\frac12\tau_b\cdot$, is a basic self-adjoint transversally elliptic differential operator, called the basic   Dirac operator.

\noindent
 We say $(M,F)$ admits a transverse \spinc structure, if $Q$ is a foliated \spinc structure, and the \spinc structure corresponds to the foliated line bundle admitting a basic connection.
Now we give a brief introduction of  basic Seiberg-Witten Theory for codimension $4$-foliation by Kordyukov, Lejmi and Weber.

\noindent
Let $(M,F)$ be a oriented closed manifold with codimension $4$ oriented Riemannian foliation $F$. Fixing a transverse $spin^c$ structure $\mathfrak s$, we consider the basic Seiberg-Witten equations, i.e.\begin{equation}
  \begin{cases}
    \Dirac^+_A\psi=0,\\
    F^+_A=q(\psi),
  \end{cases}\label{sw-equation-4}
\end{equation}
for $(A,\psi)\in \mathcal A_b(\mathfrak s)\times \Gamma_b(S^+)$ and  $q(\psi)=\psi\otimes\bar\psi-\frac12|\psi|^2$. For the second equation, we identify the    traceless endmorphism of the spinor bundle with the imaginary self-dual bundle(see Morgan' book  \cite{Morgan}). The basic gauge group is defined by
\[\mathcal G_b=\{u:M\to U(1)\big|~L_Xu\equiv0,\mbox{ for all }X\in \Gamma(F).\}\] Kordyukov, Lejmi and Weber   posed a necessary condition, i.e.   $H^1_b(M)\cap H^1(M,\mathbb Z)\subset H^1(M)$ is a lattice of $H^1_b(M)$, to show the compactness of the moduli space.
With repeat the similar argument of the classical Seiberg-Witten theory, we obtain the basic Seiberg-Witten invariant.

\noindent In this section,
let $(M,F,g^T)$ be a Riemannian foliation with a bundle-like metric. Suppose that it admits a transverse \spinc structure $\mathfrak s$.

\begin{lemma}\label{lemma-Kato}
  Let $A$ be a basic $spin^c$ connection and $\phi$ be a basic spinor field.  Then,
  $$|\nabla|\phi||\leq |\nabla^A\phi|.$$
  Let $\Dirac^A_b=\Dirac^A_{tr}-\frac12\tau_b$ be the basic Dirac operator. If moreover $\Dirac^A_b\phi=0$, $\phi\neq0$ and the above equality holds, then $\nabla^A_{tr}\phi=\frac12\kappa_b\phi$, where the right hand side denotes the scalar multiplication.
\end{lemma}
\begin{pf}
  Since $\nabla^A_{tr}$ is metric compatible, the first statement  can be obtained by applying  the arguments of {\cite[Lemma 3.1]{FZ}}. When the equality holds, one has that
  \[\nabla^A_{tr,e_i}=\alpha_i\phi, ~i=1,...,m\]
  for some real constant $\alpha_i$, where $\{e_i\}$ is the local orthonormal basis of the transverse bundle. The formula $\Dirac^A_{b}\phi=0$ implies that
  \[(\sum_i\alpha_ie^i-\frac12\kappa_b)\phi=0.\]
  Since $\phi\neq0$, we have that $\nabla^A_{tr}\phi=\frac12\kappa_b\phi$.
\end{pf}

\begin{prop}
  Let $(M,F)$ be a transverse oriented closed codimension $4$ foliation. Suppose  that $(M,F)$ admits a bundle-like metric $g$ and a transverse \spinc structure  $\mathfrak s$. If there is an irreducible solution $(A,\phi)$ to the basic Seiberg-Witten equations with respect to $g$ and $\mathfrak s$, then we get the inequality \[\lambda^T(g)\leq -\sqrt{\|F^+_A\|^2}.\]
   Moreover, the equality holds if  and only if $g$ is a taut transverse K\"ahler metric .
   \end{prop}
\begin{pf}
Set $|\phi|^2_\epsilon=|\phi|^2+\epsilon^2$,  the same argument implies that $$|\nabla|\phi|_{\epsilon}|^2\leq |\nabla^A\phi|^2.$$
 By \cite{GK}, we have the Weitzenb\"ock formula,
\begin{equation}
  (\Dirac^A_b)^2=(\nabla^A_{tr})^*\nabla^A_{tr}
  +\frac14(Scal^T+|\kappa_b|^2-2\delta_b\kappa_b)+\frac12 F_A\label{eqn-Weitzenbbock}
\end{equation}
For convenience, we set $S^T_b=Scal^T+|\kappa_b|^2-2\delta_b\kappa_b$.
Since $(A,\phi)$ is an irreducible solution to the basic Seiberg-Witten equations. By calculation, we have that
\[0=\frac12\Delta_b|\phi|^2+|\nabla^A_{tr}\phi|^2+\frac14S^T_b|\phi|^2+\frac14|\phi|^4.\]
Hence, \[\int_M|\nabla^A_{tr}\phi|^2+\frac14 S^T_b|\phi|^2=-\frac14\int_M|\phi|^4.\]
Since $\lambda^T(g^T)$ is the lowest eigenvalue of the operator $4\Delta_b+S^T_b$, one has that
\[\lambda^T(g^T)\int_M|\phi|^2_\epsilon
\leq \int_M(4|\nabla^T|\phi|_\epsilon|^2+S^T_b|\phi|^2_\epsilon)
\leq\int_M4|\nabla^A_{tr}\phi|^2+S^T_b|\phi|^2_\epsilon.\]
For $\lambda^T(g)\leq0$ and $\epsilon<<1$, by Cauchy-Schwarz inequality, it holds that
\begin{eqnarray*}
\lambda^T(g)Vol(g)^{1/2}(\int_M|\phi|^4_\epsilon)^{1/2}&\leq&
\lambda^T(g)\int_M|\phi|^2_\epsilon\\
&\leq&-\int_M|\phi|^4+\epsilon^2\int_M(S^T_b).
\end{eqnarray*}
Taking $\epsilon\to0$, we have
\[\bar\lambda^T(g)\leq-(\int_M|\phi|^4)^{1/2}.\]
Combining with the second equation in the Seiberg-Witten equations, one has that
\[\bar\lambda^T(g)\leq-(8\int_M|F^+_A|^2)^{1/2}. \]
Thus, we get
\[\bar\lambda^T(g)\leq-\sqrt{32\pi^2|c^+_1(\mathfrak s)|^2}.\]
If the equality holds, we have that $|\phi|$ is constant. Moreover, combining Lemma \ref{lemma-Kato} and the formula
\[d|\phi|^2=\kappa_b|\phi|^2=0,\]
we have that $\kappa_b=0$, i.e. the bundle-like metric is taut, and $Scal^T=-|\phi|^2:= Constant$.
 By the second  equation $F^+_A=\phi\otimes\phi^*-\frac12|\phi|^2$, we have that
$\nabla^TF^+_A=0$.
 The last equality implies that $g$ has a transverse K\"ahler structure with respect to the transverse  K\"ahler form  $\sqrt2\frac{F^+_A}{|F^+_A|}$.
\end{pf}

\noindent
Combining with the above proposition and  \cite[Claim under 2.3]{Pe1}. We have the foliated version of \cite[Theorem 1.1]{FZ}, i.e. Theorem \ref{thm-main-2}. Since the arguments of proof is no essentially different, we omit its proof here.

\vspace{5mm}

For a 5-dimensional Sasaki manifold, a fundamental topic is to determine whether this manifold admits a Sasaki-Einstein metric or not. At the end of this paper, we give a corollary for 5-dimensional Sasaki manifolds admitting a Sasaki-Einstein metric. We recall the definition of Sasaki manifold.

\begin{defi}
  A Sasakian manifold of dimension $2 n+1$ is a Riemannian manifold $\left(M^{2 n+1}, g\right)$ with the property that its metric cone $\left(C(M)=\mathbb{R}_{>0} \times M, \bar{g}=\right.$ $\left.d r^2+r^2 g\right)$ is K\"ahler. The Reeb vector field is $\xi=J\left(r \partial_r\right)$, where $J$ denotes the integrable complex structure on $C(M)$.
\end{defi}

Let $Q=TM/\xi$ be the quotient space of a Sasaki manifold $(M,g,\xi)$. It is  known that $Q$ has a transverse K\"ahler structure c.f. \cite{BG}, which induces a canonical transverse $spin^c$ structure. 
Combining with  the work of
\cite{KLW}, we have the following corollary.

\begin{cor}
  Let $(M,g,\xi)$ be a 5-dimensional closed Sasaki-manifold with $b^+_b>1$.
  Suppose that the identity $\bar\lambda^T(g)=-\sqrt{32\pi^2 c^2_b(\mathfrak s)}$ holds, where $\mathfrak s$ denotes by the canonical transverse $spin^c$ structure.   Then, $(M,g,\xi)$ is  a Sasaki-Einstein manifold.
\end{cor}



 College of Mathematics and Statistics, Chongqing University,
Huxi Campus, Chongqing, 401331, P. R. China

Chongqing Key Laboratory of Analytic Mathematics and Applications, Chongqing University, Huxi Campus, Chongqing, 401331, P. R.
China

 E-mail: lindexie@126.com


\begin{thebibliography}{99}


\bibitem{AL}
 J. A. Alvarez L\`opez, The basic component of the mean curvature of Riemannian
foliations, Ann. Global Anal. Geom.,   10(1992), 179-194.


\bibitem{BG}
C.P. Boyer and K. Galicki, Sasakian Geometry, Oxford University Press, Oxford, 2008.




\bibitem{BKR}
 J. Br\"uning, F. W. Kamber and K. Richardson, Index theory for basic Dirac operators on Riemannian foliations, Noncommutative geometry and global analysis, 39-81.
Contemp. Math., 546, Amer. Math. Soc., Providence, RI, 2011.


\bibitem{Collins}
T.C. Collins, The transverse entropy functional and the Sasaki-Ricci flow, Tran. A.M.S.,
  365(2013), no. 3,  1277-1303.

\bibitem{D}
 D. Dominguez, A tenseness theorem for Riemannian foliations, C.R. Acad.
Sci. S\'er. I, vol. 320(1995), no. 11, 1331-1335.

\bibitem{FZ}
F.Fang and Y. Zhang, Perelman's $\lambda$-functional and Seiberg-Witten equations, Front. Math. China,  vol. 2 (2007), no.2  191-210. 


\bibitem{GK}
 J. F. Glazebrook and F. W. Kamber, Transversal Dirac families in Riemannian foliations, Comm. Math. Phys., vol. 140 (1991), no. 2, 217-240.






\bibitem{Hamilton}
R. Hamilton, Three-manifolds with positive Ricci curvature, J. Diff. Geom.,  vol. 17 (1982), no. 2, 255-306.


\bibitem{KA}
A. El Kacimi-Alaoui, Operateurs transversalement elliptiques sur un feuilletage
riemannien et applications, Compositio Math., vol. 73 (1990), 571-106.

\bibitem{KT}
F.W. Kamber and P. Tondeur, Foliated bundles and characteristic classes,
Lecture Notes in Mathematics,   493, Springer-Verlag, Berlin-New York,
1975.

\bibitem{KLW}
Y. Kordyukov, M. Lejmi and P. Weber,
Seiberg-Witten invariants on manifolds with Riemannian foliations of codimension 4,
J. of Geo. and Phy., vol. 107 (2016),  114-135.









\bibitem{LMR}
 Miroslav Lovri\'c,  Min-Oo Maung and Ernst A. Ruh,  Deforming transverse Riemannian
metrics of foliations, Asian J. Math., vol.   4 (2000), no. 2, 303-314.



\bibitem{Molino}
P. Molino, Riemannian foliations, Progress in Mathematics, vol. 73, Birkh\"auser, Boston, 1988.

\bibitem{Morgan}
 J. W. Morgan, The Seiberg-Witten equations and applications to the topology of smooth
four-manifolds, Mathematical Notes 44, Princeton University Press, Princeton, NJ, 1996.


\bibitem{Pe1}
 G.Perelman, The entropy formula for the Ricci flow and its geometric applications, arXiv:math/0211159.

 \bibitem{PR}
 I. Prokhorenkov and K. Richardson,
Perturbations of basic Dirac operators on Riemannian foliations, Internat. J. Math.   24 (2013), no. 9, 1350072, 26 pp.

\bibitem{Reinhart}
 B. Reinhart, Foliated manifolds with bundle-like metrics, Ann. Math., vol.  69 (1959)(2), 119-132.


\bibitem{Rummler}
 H. Rummler, Quelques notions simples en geometrie riemannienne et leurs
applications aux feuilletages compacts, Comment. Math. Helv., vol.  54(1979), no. 2, 224-239.


\bibitem{SC}
E. Sanmartin-Carbon, The manifold of bundle-like metrics of a Riemannian foliation, Quart. J. Math. Oxford Ser. (2), vol. 48 (1997), no. 190, 243-254.

\bibitem{SWZ}
K. Smoczyk, G. Wang and Y. Zhang, The Sasaki-Ricci flow, Inter. J. of Math., vol. 21 (2010), no. 7, 951-969.

\bibitem{Ton2}
 Ph. Tondeur, Geometry of Foliations, Birkh\"auser, Basel, 1997.


\bibitem{Wang}
 S. Wang, A higher dimensional foliated Donaldson theory, I, Asian J. Math., vol. 19 (2015), no.3, 527-554.


\bibitem{W}
 E. Witten, Monopoles and four-manifolds, Math. Res. Lett., vol. 1 (1994), no. 6, 769-796.

\bibitem{Y}
S. Yorozu, Behavior of geodesics in foliated manifolds
with bundle-like metric, J. Math. Soc. Japan, vol. 35 (1983)   no. 2, 251-272.


\end{thebibliography}
\end{document}